\documentclass[letterpaper,11pt]{article}

\usepackage{amssymb,amsmath,amscd,amsfonts}
\usepackage{epsf,epsfig}
\usepackage{times,bm,mathrsfs,bbm}
\usepackage{color,graphicx}

\definecolor{Red}{rgb}{1,0,0}
\definecolor{Blue}{rgb}{0,0,1}
\definecolor{Olive}{rgb}{0.41,0.55,0.13}
\definecolor{Green}{rgb}{0,1,0}
\definecolor{MGreen}{rgb}{0,0.8,0}
\definecolor{DGreen}{rgb}{0,0.55,0}
\definecolor{Yellow}{rgb}{1,1,0}
\definecolor{Cyan}{rgb}{0,1,1}
\definecolor{Magenta}{rgb}{1,0,1}
\definecolor{Orange}{rgb}{1,.5,0}
\definecolor{Violet}{rgb}{.5,0,.5}
\definecolor{Purple}{rgb}{.75,0,.25}
\definecolor{Brown}{rgb}{.75,.5,.25}
\definecolor{Grey}{rgb}{.5,.5,.5}
\definecolor{Black}{rgb}{0,0,0}

\newtheorem{theorem}{Theorem}

\newtheorem{corollary}{Corollary}

\newtheorem{lemma}{Lemma}

\newenvironment{proof}{\noindent{\textbf{Proof:}}}{$\blacksquare$\vskip\belowdisplayskip}

\newcommand{\real}{\mathbb{R}}

\newcommand{\expec}{\mathbb{E}}

\begin{document}

\title{\vspace{-3cm} 
Reconstruction on Trees:\\ 
Exponential Moment Bounds for Linear Estimators
}
\author{
Yuval Peres\footnote{
Theory Group,
Microsoft Research.}
\and
Sebastien Roch\footnote{
Department of Mathematics, UCLA.
}
}
\maketitle

\begin{abstract}
Consider a Markov chain $(\xi_v)_{v \in V} \in [k]^V$ on the infinite
$b$-ary tree $T = (V,E)$ with irreducible edge transition
matrix $M$, where $b \geq 2$, $k \geq 2$ and $[k] = \{1,\ldots,k\}$.
We denote by $L_n$ the level-$n$ vertices of $T$.
Assume $M$ has a real second-largest (in absolute value) eigenvalue $\lambda$ 
with corresponding real eigenvector $\nu \neq 0$.
Letting $\sigma_v = \nu_{\xi_v}$, we consider the following root-state
estimator,  
which was introduced by
Mossel and Peres (2003) in the context of the ``recontruction problem''
on trees:
\begin{equation*}
S_n = (b\lambda)^{-n} \sum_{x\in L_n} \sigma_x.
\end{equation*}
As noted by Mossel and Peres, when $b\lambda^2 > 1$ (the so-called Kesten-Stigum reconstruction phase)
the quantity $S_n$ has uniformly bounded variance. Here,
we give bounds on the moment-generating
functions of $S_n$ and $S_n^2$ when $b\lambda^2 > 1$. Our results have implications for the
inference of evolutionary trees.
\end{abstract}


\section{Introduction}

We first state our main theorem. Related results and applications
are discussed at the end of the section.

\paragraph{Basic setup.}
For $b \geq 2$, let $T = (V,E)$ be the infinite $b$-ary tree rooted at $\rho$. 
Denote by $T_n$ the first $n \geq 0$ levels of $T$. 
Let $M = (M_{ij})_{i,j=1}^k$ be a $k\times k$ irreducible stochastic matrix with stationary distribution $\pi > 0$.
Assume $M$ has a real second-largest (in absolute value) eigenvalue $\lambda$ 
and let $\nu \neq 0$ be a real right eigenvector corresponding to $\lambda$
with
\begin{equation*}
\sum_{i=1}^k \pi_i \nu_i^2 = 1.
\end{equation*}
Let $[k] = \{1,\ldots,k\}$.
Consider the following Markov process on $T$:
pick a root state $\xi_\rho$ in $[k]$ according to $\pi$;
moving away from the root, apply the channel $M$ to each edge 
independently. Denote by $(\xi_v)_{v \in V}$ the state assignment
so obtained and let
\begin{equation*}
\sigma_v =  \nu_{\xi_v},
\end{equation*}
for all $v \in V$

\paragraph{Reconstruction.}
In the so-called ``reconstruction problem,'' one seeks---roughly speaking---to
infer the state at the root from the states at level $n$, as $n \to \infty$.
This problem has been studied extensively in probability theory and
statistical physics. See e.g.~\cite{EvKePeSc:00} for background and references.
Here, we are interested in the following root-state estimator introduced
in~\cite{MosselPeres:03}.
For $n \geq 0$, let $L_n$ be the vertices of $T$ at level $n$.
Consider the following quantity
\begin{equation}\label{eq:sn}
S_n = \frac{1}{(b\lambda)^n} \sum_{x\in L_n} \sigma_x.
\end{equation} 
It is easy to show that for all $n \geq 0$
\begin{equation*}
\expec[S_n\,|\,\xi_\rho] 
= \sigma_\rho,
\end{equation*}
that is, $S_n$ is ``unbiased.''
Moreover, it was shown in~\cite{MosselPeres:03} that
in the so-called Kesten-Stigum reconstruction phase, that is,
when $b \lambda^2 > 1$, it holds that for all $n \geq 0$ 
\begin{equation*}
\max_i \expec[S_n^2\,|\,\xi_\rho = i] 
\leq C  < +\infty,
\end{equation*}
where $C = C(M)$ is a constant depending only
on $M$ (not on $n$).

\paragraph{Main results.}
For $n \geq 0$, $i =1,\ldots,k$, and $\zeta \in \real$, let
\begin{equation*}
\Gamma_n^i(\zeta) = \expec[e^{\zeta S_n}\,|\,\xi_\rho = i],
\end{equation*}
and
\begin{equation*}
\widetilde\Gamma_n^i(\zeta) = \expec[e^{\zeta S^2_n}\,|\,\xi_\rho = i].
\end{equation*}
We prove the following.
\begin{theorem}[Exponential Moment Bound]\label{thm:main}
Assume $M$ is such that $b \lambda^2 > 1$.
Then, there is $c = c(M) < +\infty$ such that for all 
$n \geq 0$, $i =1,\ldots,k$, and $\zeta \in \real$,
it holds that
\begin{equation*}
\Gamma_n^i(\zeta) \leq e^{ \nu_i \zeta + c \zeta^2} < +\infty.
\end{equation*}
Note that $\nu_i = \expec[S_n\,|\,\xi_\rho = i]$.
\end{theorem}
\begin{corollary}\label{cor}
Assume $M$ is such that $b \lambda^2 > 1$.
Then, there is $\tilde\zeta = \tilde\zeta(M) \in (0, +\infty)$ 
and $\widetilde C = \widetilde C(M) < +\infty$ such that for all 
$n \geq 0$, $i =1,\ldots,k$, and $\zeta \in (-\tilde\zeta, \tilde\zeta)$,
it holds that
\begin{equation*}
\widetilde\Gamma_n^i(\zeta) \leq \widetilde C < +\infty.
\end{equation*}
\end{corollary}
The proofs of Theorem~\ref{thm:main} and Corollary~\ref{cor} can be found in Section~\ref{section:proofs}.

\paragraph{Related results.}
Moment-generating functions of random variables similar to (\ref{eq:sn}) have been studied
in the context of multi-type branching processes. In particular, Athreya and Vidyashankar~\cite{AthreyaVidyashankar:95}
have obtained large-deviation results for quantities of the type (in our setting)
\begin{equation*}
R_n = b^{-n} Z_n \cdot w - \pi \cdot w,
\end{equation*}
where $w \in \real^k$ and $Z_n = (Z_n^{(1)},\ldots,Z_n^{(k)})$ is the ``census'' vector, that is,
\begin{equation*}
Z_n^{(i)} = |\{x\in L_n\ :\ \xi_x = i\}|,
\end{equation*}
for all $i\in [k]$.
However, note that we are interested in the \emph{degenerate} case $w = \nu \perp \pi$
(see e.g.~\cite{HornJohnson:85}) and our results cannot be deduced from~\cite{AthreyaVidyashankar:95}.

Note moreover that our bounds cannot hold when $b\lambda^2 < 1$. Indeed, in that case, a classical CLT
of Kesten and Stigum~\cite{KestenStigum:66} for multi-type branching processes implies that the quantity
\begin{equation*}
Q_n \equiv (b\lambda^2)^{n/2} S_n =  \frac{1}{b^{n/2}} \sum_{x\in L_n} \sigma_x,
\end{equation*}
converges in distribution to a centered Gaussian with a finite variance (independently of the root state).
See~\cite{MosselPeres:03} for more on the Kesten-Stigum CLT and its relation to the reconstruction
problem.

\paragraph{Motivation.}
The motivation behind our results comes from mathematical biology.
More particularly, our main theorem has recently played a role in the solution of important questions 
in mathematical phylogenetics, which we now briefly discuss.

As mentioned above, the quantity $S_n$ arises naturally in the reconstruction problem 
as a simple ``linear'' estimator of the root state~\cite{EvKePeSc:00, MosselPeres:03}. 
In the past few years, deep connections have been established between the reconstruction problem
and the inference of phylogenies---a central problem in computational biology~\cite{SempleSteel:03,Felsenstein:04}.
A phylogeny is a tree representing the evolutionary history of a group of organisms, where
the leaves are modern species and the branchings correspond to past speciation events.
To reconstruct phylogenies, biologists extract (aligned) biomolecular sequences from extant species.
It is standard in evolutionary biology to model such collections of sequences as 
\emph{independent samples from the leaves of 
a Markov chain on a finite tree}
\begin{equation}\label{eq:samples}
\mathbb{S} = \{(\sigma^i_x)_{x\in L_n}\}_{i=1}^\ell,
\end{equation}
where $\ell$ is the sequence length. 
The goal of phylogenetics is to infer the \emph{leaf-labelled} tree
that generated these samples. 
In particular, developing reconstruction techniques 
that require as few samples as possible
is of practical importance. 

An insightful conjecture of Steel~\cite{Steel:01} suggests that
the reconstruction of phylogenies can be achieved from much shorter sequences when the reconstruction
problem is ``solvable,'' in particular in the Kesten-Stigum reconstruction phase. 
This conjecture has been established in the binary symmetric case (equivalent
to the ferromagnetic Ising model), that is, the case $k=2$ and $M$ symmetric, by Mossel~\cite{Mossel:04a} 
and Daskalakis et al.~\cite{DaMoRo:09}. 
The main idea behind these results is to ``boost'' standard tree-building techniques
by inferring ancestral sequences.
See~\cite{Mossel:04a,DaMoRo:09} for details. 

Establishing Steel's conjecture under more realistic models of sequence
evolution (i.e., more general transition matrices $M$) 
is a major open problem in mathematical phylogenetics. 
Roughly, to reconstruct a phylogeny from samples 
at level $n$ one iteratively joins the most correlated pairs of nodes,
starting from level $n$ and moving towards the root. To estimate the correlation
between \emph{internal} nodes $u$ and $v$ on level $m < n$  using only (\ref{eq:samples}) 
it is natural to consider quantities such as
\begin{equation}
\widehat{\mathrm{Cov}}[u, v]
= \frac{1}{\ell} \sum_{i=1}^\ell \left((b\lambda)^{-(n-m)} \sum_{x\in L^u_n} \sigma^i_x\right)
\left((b\lambda)^{-(n-m)} \sum_{x\in L^v_n} \sigma^i_x\right),
\end{equation} 
where $L_n^u$ is the set of nodes on level $n$ below $u$.
In words, we estimate the correlation between the \emph{reconstructed} states
at $u$ and $v$.
Proving concentration of such quantities necessitates uniform bounds on the moment-generating
functions of $S_n$ and $S_n^2$---our main result.
We note in particular that our main theorem was recently used by Roch~\cite{Roch:09}, building on~\cite{Roch:08},
to prove Steel's conjecture for general $k$ and reversible transition matrices
of the form $M = e^{tQ}$ in the Kesten-Stigum phase. 
Moreover, this result was established using
a surprisingly simple algorithm known in phylogenetics as a ``distance-based method,''
thereby contradicting a conjecture regarding the weakness of this widely used 
class of methods.
See~\cite{Roch:08} for background.

\paragraph{Organization.} The proof of our results can be found in Section~\ref{section:proofs}.

\section{Proof}\label{section:proofs}

We first prove our main theorem in a neighbourhood around zero.
\begin{lemma}\label{lem:main1}
Assume $M$ is such that $b \lambda^2 > 1$.
Then, there is $c' = c'(M) < +\infty$ and $\zeta_0 \in (0,+\infty)$ 
such that for all 
$n \geq 0$, $i =1,\ldots,k$, and $|\zeta| < \zeta_0$,
it holds that
\begin{equation*}
\Gamma_n^i(\zeta) \leq e^{ \nu_i \zeta + c' \zeta^2}.
\end{equation*}
\end{lemma}
\begin{proof}
We prove the result by induction on $n$.
For $n=0$, note that
\begin{equation*}
\Gamma_0^i(\zeta) = e^{\nu_i \zeta}, 
\end{equation*}
so the first step of the induction holds for all $c' > 0$
and all $\zeta \in \real$.

Now assume the result holds for $n > 0$ with $c'$ and $\zeta_0$
to be determined later.
For $n \geq 0$, $i =1,\ldots,k$, and $\zeta \in \real$, let
\begin{equation*}
\gamma_n^i(\zeta) = \ln\Gamma_n^i(\zeta).
\end{equation*}
Let $\alpha_1,\ldots,\alpha_b$ be the children of
$\rho$ and, for $\omega=1,\ldots,b$, denote by $L^\omega_{n+1}$ the descendants
of $\alpha_\omega$ on the $n+1$'st level. 
For $\omega=1,\ldots,b$, let
\begin{equation*}
S^\omega_{n+1} = \frac{1}{(b\lambda)^n}\sum_{x\in L^\omega_{n+1}} \sigma_x. 
\end{equation*}
Note that conditioned on $\xi_\rho$, the random vectors
\begin{equation*}
(\xi_x)_{x\in L^1_{n+1}}, \ldots, (\xi_x)_{x\in L^b_{n+1}},
\end{equation*}
are independent and identically distributed. Hence, the variables
\begin{equation*}
S^1_{n+1},\ldots,S^b_{n+1},
\end{equation*} 
are also conditionally independent and identically distributed.
Applying the channel to the first level of the tree and
using the induction hypothesis, we have for $\zeta \in (-\zeta_0, \zeta_0)$
\begin{eqnarray*}
\gamma_{n+1}^i(\zeta)
&=& \ln\expec[e^{\zeta S_{n+1}}\,|\,\xi_\rho = i]\\
&=& \ln\expec\left[\exp\left(\frac{\zeta}{b\lambda} \sum_{\omega=1}^b S^\omega_{n+1}\right)
\,\bigg|\,\xi_\rho = i\right]\\
&=& b\ln\expec\left[\exp\left(\frac{\zeta}{b\lambda} S^1_{n+1}\right)
\,\bigg|\,\xi_\rho = i\right]\\
&=& b\ln\left(\sum_{j=1}^k M_{i j} \expec\left[\exp\left(\frac{\zeta}{b\lambda} S^1_{n+1}\right)
\,\bigg|\,\xi_{\alpha_1} = j\right]\right)\\
&=& b\ln\left(\sum_{j=1}^k M_{i j} \Gamma^{j}_{n}\left(\frac{\zeta}{b\lambda}\right)\right)\\
&\leq& b\ln\left(\sum_{j=1}^k M_{i j} e^{\nu_{j}(\frac{\zeta}{b\lambda}) + c'(\frac{\zeta}{b\lambda})^2}\right),
\end{eqnarray*}
where we used that by assumption
\begin{equation*}
|b\lambda| \geq \frac{1}{|\lambda|} \geq 1,
\end{equation*}
so that $\zeta/(b\lambda) \in (-\zeta_0,\zeta_0)$.
By a Taylor expansion, as $\zeta_0$ goes to zero (in particular $\zeta_0 < 1$), we have 
\begin{eqnarray*}
\gamma_{n+1}^i(\zeta)
&\leq& c'\frac{\zeta^2}{b\lambda^2} \\
&& + b\ln\left(\sum_{j=1}^k M_{i j} \left[1 + \nu_{j}\left(\frac{\zeta}{b\lambda}\right) 
+ \frac{1}{2} \nu_{j}^2\left(\frac{\zeta}{b\lambda}\right)^2 
+ |\zeta|^3 \right]\right)\\ 
&\leq& c'\frac{\zeta^2}{b\lambda^2}\\
&& + b\ln\left( 1 + \lambda \nu_{i}\left(\frac{\zeta}{b\lambda}\right) 
+ \frac{1}{2} \|\nu\|_\infty^2\left(\frac{\zeta}{b\lambda}\right)^2 
+ |\zeta|^3 \right)\\
&\leq& \nu_{i}\zeta 
+ \left\{c' + \frac{1}{2} \|\nu\|_\infty^2\right\}\frac{\zeta^2}{b\lambda^2} 
- \frac{1}{2}\frac{\nu_i^2 \zeta^2}{b}
+ O_{\zeta_0}(|\zeta|^3)\\
&\leq& \nu_{i}\zeta 
+ \left\{c' + \frac{1}{2} \|\nu\|_\infty^2\right\}\frac{\zeta^2}{b\lambda^2} 
+ O_{\zeta_0}(|\zeta|^3).
\end{eqnarray*}
Choose $c' > 0$ large enough so that
\begin{equation*}
c' > \left\{c' + \frac{1}{2} \|\nu\|_\infty^2\right\}\frac{1}{b\lambda^2},
\end{equation*}
that is,
\begin{equation*}
c' > \frac{\|\nu\|_\infty^2}{2 b \lambda^2}\left(1 - \frac{1}{b\lambda^2}\right)^{-1}.
\end{equation*}
Note that $c'$ is well defined when $b\lambda^2 > 1$.
Then there is $\zeta_0 \in (0,+\infty)$ such that
for all $\zeta \in (-\zeta_0,\zeta_0)$
\begin{equation*}
\gamma_{n+1}^i(\zeta) \leq \nu_{i}\zeta + c' \zeta^2.
\end{equation*}
That concludes the proof.
\end{proof}

The following lemma deals with values of $\zeta$
away from zero.
\begin{lemma}\label{lem:main2}
Assume $M$ is such that $b \lambda^2 > 1$. Let $\zeta_0 \in (0,+\infty)$
be as in Lemma~\ref{lem:main1}.
Then, there is $c'' = c''(M) < +\infty$ 
such that for all 
$n \geq 0$, $i =1,\ldots,k$, and $|\zeta| \geq \zeta_0$,
it holds that
\begin{equation*}
\Gamma_n^i(\zeta) \leq e^{c'' \zeta^2}.
\end{equation*}
\end{lemma}
\begin{proof}
Let $c'$ be as in Lemma~\ref{lem:main1}.
Let $\zeta_1 \in (0,+\infty)$ be such that
\begin{equation}\label{eq:jump}
\zeta_1 < \frac{\zeta_0}{|b\lambda|}.
\end{equation}
Choose $c'' > c'$ large enough so that
\begin{equation}\label{eq:assump1}
e^{\nu_i \zeta + c' \zeta^2} \leq e^{c'' \zeta^2},
\end{equation}
for all $|\zeta| > \zeta_1$ and for all $i = 1,\ldots, k$.

Let $n \geq 0$ and $\zeta$ with $|\zeta| \geq \zeta_0$
be fixed. Note that, when we relate the exponential moment
at level $m$ to that at level $m-1$ with a recursion as in the proof of
Lemma~\ref{lem:main1}, the value of $\zeta$ is effectively divided
by $b\lambda$. 
Therefore,
there are two cases in the proof: either we reach
the interval $(-\zeta_0, \zeta_0)$ by the time we reach $m=0$
in the recursion; or we do not.
\begin{enumerate}
\item 
First assume that 
\begin{equation}\label{eq:assump2}
\left|\frac{\zeta}{(b\lambda)^{n}}\right| \geq \zeta_0, 
\end{equation}
that is, we do not reach $(-\zeta_0,\zeta_0)$.
We prove the result by induction on the level $m=0,\ldots,n$.
At $m=0$,
we have
\begin{equation*}
\Gamma_0^i\left(\frac{\zeta}{(b\lambda)^{n}}\right) 
= e^{\nu_i (\frac{\zeta}{(b\lambda)^{n}})} \leq e^{c'' (\frac{\zeta}{(b\lambda)^{n}})^2}, 
\end{equation*}
by (\ref{eq:assump1}) and (\ref{eq:assump2})
for all $i=1,\ldots,k$.
Assume for the sake of the induction that
\begin{equation*}
\Gamma_m^i\left(\frac{\zeta}{(b\lambda)^{n - m}}\right)
\leq e^{c'' (\frac{\zeta}{(b\lambda)^{n - m}})^2},
\end{equation*}
for all $i=1,\ldots,k$.
Using the calculations of Lemma~\ref{lem:main1}, we have
\begin{eqnarray*}
\gamma_{m+1}^i\left(\frac{\zeta}{(b\lambda)^{n - (m+1)}}\right)
&=& b\ln\left(\sum_{j=1}^k M_{i j} \Gamma^{j}_{m}
\left(\frac{1}{b\lambda}\frac{\zeta}{(b\lambda)^{n - (m+1)}}\right)\right)\\
&\leq& b\ln\left(\sum_{j=1}^k M_{i j} e^{c''(\frac{\zeta}{(b\lambda)^{n - m}})^2}\right)\\
&=& b c''\left(\frac{\zeta}{(b\lambda)^{n - m}}\right)^2\\
&=& \frac{b}{b^2 \lambda^2} c''\left(\frac{\zeta}{(b\lambda)^{n - (m+1)}}\right)^2\\
&\leq& c''\left(\frac{\zeta}{(b\lambda)^{n - (m+1)}}\right)^2,
\end{eqnarray*}
where we used $b \lambda^2 > 1$ on the last line.
The proof of the first case follows by induction, that is, we have
\begin{equation*}
\Gamma_{n}^i(\zeta) \leq e^{c'' \zeta^2},
\end{equation*}
for all $i=1,\ldots,k$.

\item
Assume now that 
\begin{equation}\label{eq:assump3}
\left|\frac{\zeta}{(b\lambda)^{n}}\right| < \zeta_0. 
\end{equation}
Let $m^*$ be the largest value in $0,\ldots,n$ such that
\begin{equation}
\left|\frac{\zeta}{(b\lambda)^{n - m^*}}\right| < \zeta_0. 
\end{equation}
The purpose of Assumption (\ref{eq:jump}) above is to make sure that we never ``jump''
entirely over the subset of $(-\zeta_0,\zeta_0)$ where
(\ref{eq:assump1}) holds. Indeed,
by (\ref{eq:jump}) and
\begin{equation}
\left|\frac{\zeta}{(b\lambda)^{n - (m^*+1)}}\right| \geq \zeta_0,
\end{equation}
it follows that we must also have
\begin{equation}
\left|\frac{\zeta}{(b\lambda)^{n - m^*}}\right| > \zeta_1. 
\end{equation}
Hence, by (\ref{eq:assump1}) and Lemma~\ref{lem:main1},
we get
\begin{equation*}
\Gamma_{m^*}^i\left(\frac{\zeta}{(b \lambda)^{n - m^*}}\right) \leq e^{c'' (\frac{\zeta}{(b \lambda)^{n - m^*}})^2},
\end{equation*}
for all $i=1,\ldots,k$.
The proof then follows by induction as in the first case above.
\end{enumerate}
\end{proof}

\noindent{\textbf{Proof of Theorem~\ref{thm:main}:}}
Let $\zeta_0$, $c'$ and $c''$ be as in Lemmas~\ref{lem:main1}
and~\ref{lem:main2}. Choose $c > c'' (> c')$ large enough so that
\begin{equation}
e^{c'' \zeta^2} \leq e^{\nu_i \zeta + c \zeta^2},
\end{equation}
for all $|\zeta| \geq \zeta_0$ and for all $i = 1,\ldots, k$.
The result then follows by combining Lemmas~\ref{lem:main1}
and~\ref{lem:main2}.
$\blacksquare$\vskip\belowdisplayskip

\noindent{\textbf{Proof of Corollary~\ref{cor}:}}
We use a standard trick relating the exponential moment of the square
to that of a Gaussian.
Let $X$ be a standard normal. Using Theorem~\ref{thm:main}
and applying Fubini we have for all $n \geq 0$
and $i=1,\ldots,k$
\begin{eqnarray*}
\expec[e^{\zeta S^2_n}\,|\,\xi_\rho = i]
&=& \expec[e^{\sqrt{2\zeta}  S_n X}\,|\,\xi_\rho = i]\\
&\leq& \expec[e^{\nu_i \sqrt{2\zeta}X + c 2 \zeta X^2}\,|\,\xi_\rho = i].
\end{eqnarray*}
The last expectation is finite for $\zeta$ small enough.
$\blacksquare$\vskip\belowdisplayskip

\bibliographystyle{alpha}
\bibliography{full}

\begin{thebibliography}{EKPS00}

\bibitem[AV95]{AthreyaVidyashankar:95}
K.~B. Athreya and A.~N. Vidyashankar.
\newblock Large deviation rates for branching processes. {II}. {T}he multitype
  case.
\newblock {\em Ann. Appl. Probab.}, 5(2):566--576, 1995.

\bibitem[DMR09]{DaMoRo:09}
Constantinos Daskalakis, Elchanan Mossel, and S{\'e}bastien Roch.
\newblock Evolutionaty trees and the {I}sing model on the {B}ethe lattice: a
  proof of {S}teel's conjecture.
\newblock Preprint, 2009.

\bibitem[EKPS00]{EvKePeSc:00}
W.~S. Evans, C.~Kenyon, Y.~Peres, and L.~J. Schulman.
\newblock Broadcasting on trees and the {I}sing model.
\newblock {\em Ann. Appl. Probab.}, 10(2):410--433, 2000.

\bibitem[Fel04]{Felsenstein:04}
J.~Felsenstein.
\newblock {\em Inferring Phylogenies}.
\newblock Sinauer, New York, New York, 2004.

\bibitem[HJ85]{HornJohnson:85}
Roger~A. Horn and Charles~R. Johnson.
\newblock {\em Matrix analysis}.
\newblock Cambridge University Press, Cambridge, 1985.

\bibitem[KS66]{KestenStigum:66}
H.~Kesten and B.~P. Stigum.
\newblock Additional limit theorems for indecomposable multidimensional
  {G}alton-{W}atson processes.
\newblock {\em Ann. Math. Statist.}, 37:1463--1481, 1966.

\bibitem[Mos04]{Mossel:04a}
E.~Mossel.
\newblock Phase transitions in phylogeny.
\newblock {\em Trans. Amer. Math. Soc.}, 356(6):2379--2404, 2004.

\bibitem[MP03]{MosselPeres:03}
E.~Mossel and Y.~Peres.
\newblock Information flow on trees.
\newblock {\em Ann. Appl. Probab.}, 13(3):817--844, 2003.

\bibitem[Roc08]{Roch:08}
S{\'e}bastien Roch.
\newblock Sequence-length requirement for distance-based phylogeny
  reconstruction: Breaking the polynomial barrier.
\newblock In {\em FOCS}, pages 729--738, 2008.

\bibitem[Roc09]{Roch:09}
S{\'e}bastien Roch.
\newblock Phase transition in distance-based phylogeny reconstruction.
\newblock Preprint, 2009.

\bibitem[SS03]{SempleSteel:03}
C.~Semple and M.~Steel.
\newblock {\em Phylogenetics}, volume~22 of {\em Mathematics and its
  Applications series}.
\newblock Oxford University Press, 2003.

\bibitem[Ste01]{Steel:01}
M.~Steel.
\newblock {My Favourite Conjecture}.
\newblock Preprint, 2001.

\end{thebibliography}

\end{document}